\documentclass[11pt]{article}
\usepackage{IEEEtrantools}
\usepackage{mathtools, nccmath}
\usepackage{latexsym}
\usepackage{amsfonts}
\usepackage{amssymb}
\usepackage{amsmath,amsfonts,amsthm}
\usepackage{mathrsfs}
\usepackage{enumerate}
\usepackage{dsfont}
\usepackage{hyperref}
\usepackage{graphicx}
\usepackage{epstopdf}
\usepackage{tikz}
\usetikzlibrary{shapes,arrows.meta,positioning}
\usepackage{cleveref}
\usepackage{tikz-cd}
\usepackage{psfrag}

\usepackage{multirow}
\usepackage{tikz}
\usetikzlibrary{arrows,decorations.markings}
\usepackage{bookmark}
\usepackage{hyperref}
\hypersetup{
     colorlinks   = true,
     citecolor    = blue
}

\newcommand{\newsection}[1]{\setcounter{equation}{0}
\setcounter{dfn}{0}
\section{#1}}

\newtheorem{dfn}{Definition}[section]
\newtheorem{thm}[dfn]{Theorem}
\newtheorem{lmma}[dfn]{Lemma}
\newtheorem{ppsn}[dfn]{Proposition}
\newtheorem{crlre}[dfn]{Corollary}

\newtheorem{rmrk}[dfn]{Remark}

\newcommand{\bbc}{\mathbb{C}}

\newcommand{\bbn}{\mathbb{N}}

\def \qed { \mbox{}\hfill
$\Box$\vspace{1ex}}

\usepackage{umoline}
\usepackage{mathtools}

\title{Sections and Chapters}

\begin{document}

\author{Keshab Chandra Bakshi and Satyajit Guin}

\title{On Pimsner-Popa orthonormal basis and  Popa's relative dimension of projections}
\maketitle

%%%%%%%%%%%%%%%%%%%%%%%%%%%%%%%%%%
%%%%%%%%%%%   ABSTRACT    %%%%%%%%%%%%%%%%
%%%%%%%%%%%%%%%%%%%%%%%%%%%%%%%%%%

\begin{abstract}
We show that any depth 2 subfactor with a simple first relative commutant has a unitary orthonormal basis. As a pleasant consequence, we produce new elements in the set of Popa's relative dimension of projections for such subfactors. We also construct infinitely many new elements in the set of relative dimension of projections for subfactors arising from complex Hadamard matrices and bi-unitary matrices.
\end{abstract}
\bigskip

{\bf AMS Subject Classification No.:} {\large 46}L{\large 37}, {\large 46}L{\large 10}.
\smallskip

projections, Jones index, commuting square

\hypersetup{linkcolor=blue}

%\tableofcontents

%%%%%%%%%%%%%%%%%%%%%%%%%%%%%%%%%%%%%%%
%%%%%%%%%%%%%%%%%%%%%%%%%%%%%%%%%%%%%%%
%%%%%%%%%%%%%%%%%%%%%%%%%%%%%%%%%%%%%%%

\newsection{Introduction}

For a subfactor $N\subset M$ of type $II_1$ factors and a projection $p\in M$, the element $E_N(p)$ may be regarded as the dimension of the projection $p$ relative to the subfactor $N$. The interesting case is when $E_N(p)$ is a scalar multiple of the identity. Indeed, motivated by the orthogonalization problem, in \cite{popa} Popa initiated the study of the set $\Lambda(M,N):=\{\alpha\in[0,1]:E_N(p)=\alpha 1,\,p\in\mathcal{P}(M)\}$. This set is an invariant for the inclusion $N\subset M$ that is closely related to the index $[M:N]$. More precisely, $\inf\big(\Lambda(M,N)\setminus\{0\}\big)=[M:N]^{-1}$. Popa computed $\Lambda(M,N)$ when $[M:N]\leq 4$ and $N$ is locally trivial; and some parts of $\Lambda(M,N)$ when $[M:N]>4$. In all these cases, it turns out that $\Lambda$ depends only on the index $[M:N]$. To be more precise regarding the case when $[M:N]>4$, let $0<t<\frac{1}{2}$ be the real number given by the quadratic equation $t(1-t)=[M:N]^{-1}$. Then, the set $\Lambda(M,N)\cap(t,1-t)$ is mysterious. In fact, Popa has asked the following question in \cite{popa} (see also \cite{pop2}):
\smallskip

\noindent\textbf{Question A:}(\cite{popa, pop2}) Is it true that if $\,N^\prime\cap M=\bbc\mbox{ and }[M : N]>4$, then $\Lambda(M,N)\cap(t,1-t)\neq\emptyset$, where $0<t<\frac{1}{2}$ is the real number given by the quadratic equation $t(1-t)=[M:N]^{-1}$?
\smallskip
 
Answer to the above question in its full generality seems very difficult. In this article, we show that the question has affirmative answer in the following situations.

\begin{thm}
Let $N\subset M$ be a finite index subfactor such that $[M:N]>4$, and $0<t<\frac{1}{2}$ be the real number given by the quadratic equation $t(1-t)=[M:N]^{-1}$. Then, the set $\Lambda(M,N)\cap(t,1-t)$ is non-empty in the following cases\,:
\begin{enumerate}[$(i)$]
\item if $N\subset M$ is a depth $2$ subfactor such that $N^\prime\cap M$ is simple;
\item if $N\subset M$ is a spin model subfactor;
\item if $N\subset M$ is a vertex model subfactor.
\end{enumerate}
\end{thm}

Moreover, we prove the following stronger version of the above theorem in some situations.
\begin{thm}
\begin{enumerate}[$(i)$]
\item Let $N\subset M$ be an even index, say $2n$ with $n\geq 3$, spin model subfactor. There are increasing sequences $\{\alpha_m^{(i)}\}_{m\geq 1}$, where $0\leq i\leq 2n-4$, in $\Lambda(M,N)\cap(t,1-t)$ with limit point $\frac{2+i}{[M:N]}$, where $t(1-t)=[M:N]^{-1}$.
\item Let $N\subset M$ be a vertex model subfactor of index $(2n)^2,\,n\geq 3$. There are increasing sequences $\{\beta_m^{(i)}\}_{m\geq 1}$, where $0\leq i\leq 2n-4$, in $\Lambda(M,N)\cap(t,1-t)$ with limit point $\frac{2+i}{\sqrt{[M:N]}}$, where $t(1-t)=[M:N]^{-1}$.
\end{enumerate}
\end{thm}

Popa has shown that if the inclusion $N\subset M$ splits $R$ and if $\Lambda(M,N)\cap(t,1-t)\neq\emptyset$, then the cardinality of this set is at least $\aleph$ \cite{popa}. This is achieved by constructing recursively a decreasing sequence with limit point $t$. Therefore, applying Popa's construction to each element of the sequence in the above theorem, we obtain lots of elements in the desired set for the cases of spin/vertex model subfactors. However, the sequence that we have constructed consists entirely of rational numbers. This leads us to propose the following relatively easier version of Popa's question\,:
\smallskip

\noindent\textbf{Question B:} For an irreducible subfactor $N\subset M$ with $[M : N]>4$, can $\Lambda(M,N)\cap(t,1-t)$ (where $t(1-t)=[M:N]^{-1}$) contain irrational numbers?
\smallskip

The difficulty in answering the above question lies in the fact that it is not known whether the set $\Lambda(M,N)$ is closed, and hence we fail to apply any limit argument. Also, it is not known whether $\Lambda$ is uncountable, or if it can contain intervals. This seems to a be a challenging problem.

To compute $\Lambda$, we have used Pimsner-Popa orthonormal basis as a technical tool. Indeed, restrictions on the maximal number of unitaries that can appear in an orthonormal basis of $M/N$ can be used to obtain a restriction to the possible values in $\Lambda(M,N)$, as observed in \cite{pop2}. However, the following question due to Popa seems extremely hard.
\smallskip

\noindent\textbf{Question C:}(\cite{pop2}) For an irreducible subfactor $N\subset M$ with index $n\geq 5$, does there exist an orthonormal basis consisting of $n$ many unitary elements?
\smallskip

For a recent progress on the above question, see \cite{BG, CKP}. As a small contribution to the above question of Popa, in this paper we have proved the following.

\begin{thm}\label{to prove}
A depth $2$ subfactor with simple first relative commutant has a unitary orthonormal basis. 
\end{thm}

%%%%%%%%%%%%%%%%%%%%%%%%%%%%%%%%%%%%%%%
%%%%%%%%%%%%%%%%%%%%%%%%%%%%%%%%%%%%%%%
%%%%%%%%%%%%%%%%%%%%%%%%%%%%%%%%%%%%%%%

\newsection{Preliminaries}\label{preliminaries}

Generalizing the classical notion of independence for $\sigma$-algebras, Popa introduced a notion of orthogonality for pairs of von Neumann subalgebras $\mathcal{P}$ and $\mathcal{Q}$ of a finite von Neumann algebra $\mathcal{M}$ \cite{P}. As a slight generalization of the notion of orthogonality, Popa also introduced the notion of `commuting square', which turns out to be an indispensable tool in subfactor theory. Consider an inclusion of finite von Neumann algebras $\mathcal{N}\subset \mathcal{M}$ with a fixed trace $tr$ on $\mathcal{M}$, and intermediate von Neumann subalgebras $\mathcal{P}$ and $\mathcal{Q}$. Thus, we obtain a quadruple of von Neumann algebras
\[
\begin{matrix}
\mathcal{Q} &\subset & \mathcal{M} \cr \cup &\ &\cup\cr \mathcal{N} &\subset & \mathcal{P}\,.
\end{matrix}
\]
If $\mathcal{P}\vee \mathcal{Q}=\mathcal{M}$ and $\mathcal{P}\wedge \mathcal{Q}=\mathcal{N}$, then a quadruple is called a {\it quadrilateral}.

\begin{dfn}[\cite{Po2, GHJ, JS}]\label{def of comm}
A quadruple
\[
\begin{matrix}
\mathcal{Q} &\subset & \mathcal{M} \cr \cup &\ &\cup\cr \mathcal{N} &\subset & \mathcal{P}
\end{matrix}
\]
of finite von Neumann algebras is called a commuting square if $E^{\mathcal{M}}_{\mathcal{P}}E^{\mathcal{M}}_{\mathcal{Q}}=E^{\mathcal{M}}_{\mathcal{Q}}E^{\mathcal{M}}_{\mathcal{P}}=E^{\mathcal{M}}_{\mathcal{N}}$. The quadruple is said to be non-degenrate if ${\overline{\mathcal{PQ}}}^{\,\tiny \mbox{SOT}}={\overline{\mathcal{QP}}}^{\,\tiny \mbox{SOT}}=\mathcal{M}$. A quadruple is called a non-degenerate commuting square (or symmetric commuting square) if it is a commuting square and non-degenerate.
\end{dfn}

For brevity, we shall sometimes write $(\mathcal{N}\subset\mathcal{P,Q}\subset \mathcal{M})$ to mean the quadruple
\[
\begin{matrix}
\mathcal{Q} &\subset & \mathcal{M} \cr \cup &\ &\cup\cr \mathcal{N} &\subset & \mathcal{P}
\end{matrix}
\]
The following result is extremely useful.
\begin{thm}[\cite{Po2}]\label{combasis}
Suppose we have a quadruple $(\mathcal{N}\subset\mathcal{P,Q}\subset \mathcal{M})$ of finite von Neumann algebras with a fixed trace $tr$ on $\mathcal{M}$. If $E^{\mathcal{M}}_{\mathcal{Q}}$ and $E^{\mathcal{P}}_{\mathcal{N}}$ are the $tr$-preserving conditional expectations, then an orthonormal Pimsner-Popa basis for $\mathcal{P}/\mathcal{N}$ via $E^{\mathcal{P}}_{\mathcal{N}}$ is also an orthonormal Pimsner-Popa basis for $\mathcal{M}/\mathcal{Q}$ via $E^{\mathcal{M}}_{\mathcal{Q}}$.
\end{thm}

A rich source of commuting squares are the complex Hadamard matrices and bi-unitary matrices. In this article, we have a special focus on two classes of subfactors, namely the spin model (also known as Hadamard subfactors) and vertex model subfactors. We briefly recall the construction of these classes of subfactors from \cite{JS}. These subfactors arise as iterated basic construction of non-degenerate commuting squares of finite-dimensional algebras. Nothing much is known about these subfactors, and these are known to be difficult. Recently, an in-depth study for `pairs' of spin model subfactors is carried out in \cite{KS1, KS2}, and several interesting results have been obtained.
\smallskip

\noindent\textbf{Notation $(1)$:~} Throughout the article, we denote $M_n(\bbc)$ simply by $M_n$ for brevity.
\smallskip

Let $u$ be a complex Hadamard matrix of order $n\times n$. Then, we have the following non-degenerate commuting square
\[
\begin{matrix}
\mbox{Ad}_u(\Delta_n) &\subset & M_n\\
\cup & &\cup\\
\bbc &\subset & \Delta_n
\end{matrix}
\]
where $\Delta_n\subset M_n$ is the algebra of diagonal matrices. Iterating the basic construction vertically, we obtain the spin model subfactor $R_u\subset R$ from the following grid of finite-dimensional inclusions
\begin{IEEEeqnarray*}{lCl}
\begin{matrix}
\Delta_n &\subset & M_n &\subset^{\,e_2} & \Delta_n\otimes M_n &\subset^{\,e_3} & M_n\otimes M_n &\subset & \cdots\cdots\cr
\cup & & \cup & & \cup & & \cup & & \cr
\mathbb{C} &\subset &  u\Delta_n u^* &\subset^{e_2} & \langle u\Delta_nu^*,e_2\rangle &\subset^{e_3} & \langle \langle u\Delta_nu^*,e_2\rangle,e_3\rangle &\subset & \cdots\cdots
\end{matrix}
\end{IEEEeqnarray*}
where $e_j$'s are the Jones projections, and the notation $\langle\mathcal{A},x\rangle$ for any inclusion of algebras $\mathcal{A\subset B}$ and $x\in\mathcal{B}$ denotes the subalgebra of $\mathcal{B}$ generated by $\mathcal{A}$ and $x$. Note that $[R:R_u]=n$. See \cite{JS} for details. Importance of this class of subfactors has been emphasized by Jones; however, not much is known about this class of subfactors. The subfactor $R_u\subset R$ is irreducible and index is $n$. It is also known that the second relative commutant is always abelian.

Now, let $v$ be a bi-unitary matrix in $M_n\otimes M_k$. Recall that a unitary matrix $w=(w_{\alpha a}^{\beta b}) \in M_n \otimes M_k$ is called a bi-unitary matrix if its block transpose $\widetilde{w}= (\widetilde{w}_{\alpha a}^{\beta b})$, with $\widetilde{w}_{\alpha a}^{\beta b}:=w_{\beta a}^{\alpha b}$, is also a unitary matrix in $M_n \otimes M_k$. In this case, we have the following nondegenerate commuting square
\[
\begin{matrix}
\mbox{Ad}_v(M_n\otimes\bbc) &\subset & M_n\otimes M_k\\
\cup & &\cup\\
\bbc &\subset & \bbc\otimes M_k
\end{matrix}
\]
and iterating the basic construction vertically, as in the case of spin model, we obtain the vertex model subfactor $R_v\subset R$ such that $[R:R_v]=k^2$. See \cite{JS} for details. Note that vertex model subfactors need not be irreducible, and its index is always square of an integer.

%%%%%%%%%%%%%%%%%%%%%%%%%%%%%%%%%%%%%%%
%%%%%%%%%%%%%%%%%%%%%%%%%%%%%%%%%%%%%%%
%%%%%%%%%%%%%%%%%%%%%%%%%%%%%%%%%%%%%%%

\newsection{A quick overview on Pimsner-Popa orthonormal basis}\label{Sec 1}

In this paper, we are interested in subfactors $N\subset M$ such that $[M:N]<\infty$. We have the unique $tr$-preserving conditional expectation from $M$ onto $N$, denoted by $E_N$. Suppose that $e_1$ is the Jones projection and $M_1$ denotes the corresponding basic construction. We begin by recalling a few necessary facts to be crucially used in this article. 

\begin{dfn}\label{1}
Let $\mathcal{N}\subset \mathcal{M}$ be a unital inclusion of von Neumann algebras equipped with a faithful normal conditional expectation $\mathcal{E}$ from $\mathcal{M}$ onto $\mathcal{N}$. Then, a finite set $\mathcal{B}:=\{\lambda_1,\ldots,\lambda_n\}\subset\mathcal{M}$ is called a  (\textit{right) Pimsner-Popa basis} for $\mathcal{M}$ over $\mathcal{N}$ via $\mathcal{E}$ if every $x\in \mathcal{M}$ can be expressed as $x=\sum_{j=1}^n\lambda_j\mathcal{E}(\lambda^*_jx)$. Equivalently (by taking adjoints), every $x\in \mathcal{M}$ can be expressed as $x=\sum_{j=1}^n \mathcal{E}(x \lambda_j)\lambda^*_j.$ Further, such a basis $\{\lambda_i\}$ is said to be {\it orthonormal} if $\mathcal{E}(\lambda^*_i \lambda_j)
= \delta_{i,j}$ for all $i, j$. A finite set $\{\lambda_1,\cdots,\lambda_n\}\subset M$ is called a left Pimsner-Popa basis if $\{\lambda^*_1,\cdots,\lambda^*_n\}\subset M$ is a right Pimsner-Popa basis. Also, a left basis $\{\lambda_i\}$ is said to be orthonormal Pimsner-Popa basis if $\mathcal{E}(\lambda_i \lambda^*_j)=\delta_{i,j}$ for all $i, j$. A right orthonormal basis $\{\lambda_i\}$ is said to be a two-sided orthonormal basis if it is simultaneously a left orthonormal basis.
\end{dfn}

\noindent\textbf{Notation $(2)$:} Unless specified otherwise, in this article by a basis for $M/N$ we mean a right Pimsner-Popa orthonormal basis for $M/N$ via $E_N$.
\smallskip
 
An important result by Pimsner and Popa states that given a finite index subfactor $N\subset M$, there always exists an orthonormal Pimsner-Popa basis for $M/N$ with respect to the unique $tr_M$-preserving conditional expectation $E_N$. The following is a well-known open problem in subfactor theory. For more details, the reader is referred to \cite{BG2}. 
\smallskip
 
\noindent\textbf{Question D:}(Jones) Does an extremal finite index subfactor always have a two-sided orthonormal basis?
\smallskip
 
This is related to Question C asked by Popa (see \cite{pop2}), mentioned in the introduction.
\smallskip

\noindent\textbf{Question E:} Is it true that an integer index subfactor has a unitary orthonormal basis if and only if it has a two-sided basis?
\smallskip

In \cite{BG}, it is observed  that if $M/N$ has a unitary basis, then $M_1/M$ will have a two-sided basis. In the converse direction, we can say the following.
\begin{thm}
Suppose $N\subset M$ is an integer index subfactor. If $M/N$ has a two-sided orthonormal basis, then $M_1/N$ has a unitary orthonormal basis.
\end{thm}
\begin{prf}
Suppose $[M:N]=n$ and let $\{\lambda_i:0\leq i\leq n-1\}$ be a two-sided orthonormal basis for $M/N$. Thus, $E^M_N(\lambda_u\lambda^*_u)=E^M_N(\lambda^*_u\lambda_u)=1$.
For $0\leq j\leq n-1$, we put $\mu_j=\sum_{k=0}^{n-1} \omega^{jk} \lambda_ke_1\lambda^*_k$. Then, for $0\leq s,t\leq n-1$ we have
\begin{align*}
\mu^*_s\mu_t & =\sum_{u,v=0}^{n-1}\omega^{-su}\omega^{tv}\lambda_ue_1\lambda^*_u\lambda_ve_1\lambda^*_v\\
&= \sum_{u,v=0}^{n-1}\omega^{-su}\omega^{tv}\lambda_u E_N(\lambda^*_u\lambda_v)e_1\lambda^*_v\\
 &= \sum_{u=0}^{n-1}\omega^{(t-s)u}\lambda_ue_1\lambda^*_u
\end{align*}
Since, $\sum_u\lambda_ue_1\lambda^*_u=1$, we conclude that each $\mu_s$ is a unitary element in $M_1$. On the other hand, for $s\neq t,$ as $E^M_N(\lambda_u\lambda^*_u)=1$, we see that 
$$E^{M_1}_N(\mu^*_s\mu_t)= {[M:N]}^{-1} \sum_{u=o}^{n-1} \omega^{(t-s)u} E^{M_1}_N(\lambda_u\lambda^*_u)= {[M:N]}^{-1} \sum_{u=0}^{n-1}\omega^{(t-s)u}=0.$$
Therefore, $\{\mu_j:0\leq j\leq n-1\}$ is a unitary orthonormal basis for $M_1/N$.\qed
\end{prf}
\smallskip

\noindent \textbf{Recent progress:} Existence of a two-sided basis for an extremal finite index subfactor can be thought of as an analogue of Hall's Marriage theorem (see \cite{BG2}). Recall that a sufactor $N\subset M$ is said be regular if the normalizer $\mathcal{N}_M(N)$ generates the von Neumann algebra $M$. In \cite{BG2}, we have proved the following.

\begin{thm}[\cite{BG2}]\label{bakved}
Every finite index regular subfactor has a two-sided orthonormal basis.
 \end{thm}
 
In \cite{BG2}, we have proved that for a finite index subfactor $N\subset M$, the generalized Weyl group $G=\frac{\mathcal{N}_M(N)}{\mathcal{U}(N)\mathcal{U}(N^{\prime}\cap M)}$ is finite. \Cref{bakved} has the following pleasant consequence. 

\begin{thm}[\cite{BG2}]
If $N\subset M$ is a finite index regular subfactor, then $$[M:N]=|G|\,\mathrm{dim}_{\mathbb{C}}(N^{\prime}\cap M).$$
\end{thm}
 
Given a subfactor $N\subset M$ with $[M:N]<\infty,$ the corresponding relative commutant $N^{\prime}\cap M$ is a finite-dimensional $C^*$-algebra. Iterating the basic construction, we obtain a tower of $II_1$ factors $N\subset M\subset M_1\subset \cdots$. In particular, we have a subfactor inclusion $N\subset M_k$ with $[M_k:N]<\infty$, and therefore, the relative commutants $N^{\prime}\cap M_k$ are also finite-dimensional $C^*$-algebras. In other words, we have a natural inclusion of finite-dimensional $C^*$-algebra $N^{\prime}\cap M_k \subset N^{\prime}\cap M_{k+1}$ for each $k$. Recall that a subfactor is said be of depth `$k$' if $k$ is the least integer such that $N'\cap M_{k-2} \subset N'\cap M_{k-1} \subset N'\cap M_k$ is an instance of basic construction of finite-dimensional $C^*$-algebras.  
 
As a follow-up to \cite{BG2}, we have proved in \cite{BG} the following:
\begin{thm}[\cite{BG}]
Any depth $2$ subfactor with simple relative commutant has a two-sided orthonormal Pimsner-Popa basis. 
\end{thm}  
On the other hand, the importance of the existence of a unitary orthonormal basis for an integer index subfactor has been outlined in \cite{popa2}. Recently, in \cite{BG} it has been proved that any finite index regular subfactor with either simple or abelian relative commutant must have a unitary orthonormal basis. Also, in \cite{BG} it is conjectured that any regular subfactor will have a unitary orthonormal basis. Indeed, in \cite{CKP} this conjecture has been verified by proving the following.
  
\begin{thm}[\cite{CKP}]
Any regular subfactor with finite Jones index has a unitary orthonormal basis. 
\end{thm}
As pointed out in \cite{BG}, this enables us to conclude that any finite index regular subfactor has depth at most $2$. For more applications of unitary basis problem, see \cite{CCKL}.
\smallskip

We conclude this section mentioning a related invariant that is recently introduced by Popa while discussing paving size of a subfactor \cite{popa2}. This is  a new invariant for a subfactor $N\subset M$, denoted by $d_{ob}(N\subset M)$, defined as the infimum of $\lVert \sum_j m^*_jm_j\rVert$ over all orthonormal basis $\{m_j\}$ for $M$ over $N$. In \cite{popa2}, Popa have shown that
\[
[M:N]\leq d_{ob}(N\subset M)\leq [M:N](\lceil[M:N]\rceil-1).
\]

This is closely related to the two-sided basis problem. Indeed, if a subfactor $N\subset M$ has either a two-sided basis or a unitary orthonormal basis, then $d_{ob}(N\subset M)=[M:N]$.
\smallskip

\noindent\textbf{Question F:} Is there any extremal finite index subfactor $N\subset M$ such that $d_{ob}(N\subset M)\neq [M:N]$?
\smallskip

We observe that if $N\subset M$ is either a regular subfactor or a depth 2 subfactor with either simple or abelian relative commutant, then $d_{ob}(N\subset M)=[M:N]$.

%%%%%%%%%%%%%%%%%%%%%%%%%%%%%%%%%%%%%%%
%%%%%%%%%%%%%%%%%%%%%%%%%%%%%%%%%%%%%%%

\newsection{Depth 2 subfactors and the existence of unitary basis}\label{Sec 2}

The goal of this section is to prove \Cref{to prove}. This makes small progress on Popa's question on existence of a unitary orthonormal basis stated earlier.

Examples of depth 2 subfactors with simple first relative commutant can be easily constructed. For example, let $N\subset M$ be an irreducible depth 2 subfactor. Then, $\mathbb{C}\otimes N\subset M_n(\mathbb{C})\otimes M$ is a depth 2 subfactor with simple relative commutant. To prove \Cref{to prove}, we need the following well-known important result.

\begin{ppsn}[Section 4.6, \cite{GHJ}]\label{wellknown}
Let $N \subset M$ be a finite-index subfactor with depth $k$. Then, $tr_M|_{N^{\prime}\cap M_{l-1}}$ is a Markov trace for the inclusion $N'\cap M_{l-2}\subset N'\cap M_{l-1}$ with modulus $[M:N]$ for all $l \geq k$.
\end{ppsn}

Before presenting the proof of \Cref{to prove}, let us briefly pause to discuss the unitary basis problem for a `connected inclusion' of finite-dimensional $C^*$-algebras which will be used crucially in the proof. The grid of relative  commutants often encodes rich structure of the corresponding subfactor.  In view of the fact that for a subfactor $N\subset M$ with  $[M:N]<\infty$, the relative commutants $N^{\prime}\cap M_k \subset N^{\prime}\cap M_{k+1}$ are connected inclusions of finite-dimensional $C^*$-algebras, the unitary basis problem in the finite-dimension seems to be of interest (see \cite{Bak}). This can be thought of as a generalization of unitary error basis in quantum information theory (see \cite{BG, CKP, BB}).

\begin{ppsn}[Propn. $3.2.3$, \cite{JS}]\label{markov-equivalence}
Let $B \subset A$ be a unital inclusion of finite dimensional $C^*$-algebras with inclusion matrix $\Lambda$. Let $\tau$ be a tracial state on $A$, $A_1:=\langle A, e_1\rangle \subset B(L^2(A, \tau))$ denote the basic construction for $B \subset A$ w.r.t. $\tau$, $E_1: A_1 \to A$ denote the unique conditional expectation induced by $\tau $ and $\beta$ be a positive real number. Then, the following statements are equivalent:
\begin{enumerate}[$(i)$]
\item There exists a tracial state $\tau_1$ on $A_1$ such that
$\tau_1(x) = \tau(x)$ for all $x\in A$ and
\[
\beta \tau_1(xe_1) = \tau(x) 
\]
for all $x \in A$.
\item There exists a tracial state $\tau_1$ on $A_1$ such that $\tau_1(x) = \tau(x)$ for all $x\in A$ and
\[
E_1(e_1) = \beta^{-1} 1.
\]
\item If $\bar{t}$ is the trace vector for $\tau$, then $\Lambda^t\Lambda \bar{t} = \beta \bar{t}$.
\end{enumerate}
\end{ppsn}

\begin{dfn}[\cite{ GHJ, JS}]
Let $B \subset A$ be a unital inclusion of finite dimensional $C^*$-algebras with inclusion matrix $\Lambda$. Let $\tau$ be a tracial state on $A$ and $A_1:=\langle A, e_1\rangle \subset B(L^2(A, \tau))$ denote the basic construction for $B \subset A$ w.r.t. $\tau$, $E_1: A_1 \to A$ denote the unique conditional expectation induced by $\tau $ and $\beta$ be a positive real number. Then, $\tau$ is said to be a Markov trace for the inclusion $B \subset A$ with modulus $\beta$ if any of equivalent conditions of \Cref{markov-equivalence} holds.
 \end{dfn}
Therefore, for a unital inclusion $B \subset A$ of finite dimensional $C^*$-algebras with inclusion matrix $\Lambda$, a tracial state $tr$ on $A$ is said to be a Markov trace for $B \subset A$ if
\[
\Lambda^t \Lambda  \bar{t} = \|\Lambda\|^2 \bar{t},
\]
where $\bar{t}$ denotes the trace vector of the tracial state $tr$. For more on Markov trace, see \cite{JS}. Recall that if $B \subset A$ is a connected inclusion (that is, the Bratelli diagram is connected), then there is a unique Markov trace for $B \subset A$. 

As a tool to prove the existence of a unitary orthonormal basis for a finite index regular subfactor, we have the following result.

\begin{ppsn}[\cite{BG}]\label{unitary-onb}
Let $A$ be a finite dimensional $C^*$-algebra which is either simple or commutative. Then, $A/\mathbb{C}$ has a unitary orthonormal basis with respect to the Markov trace for the unital inclusion $\mathbb{C} \subset A$.
\end{ppsn}

Recently in \cite{CKP}, the above has been generalized to the following.

\begin{thm}[\cite{CKP}]\label{kribs}
Let $A$ be any finite-dimesnional $C^*$-algebra. Then, $A/\mathbb{C}$ has a unitary orthonormal basis with respect to the Markov trace for $\mathbb{C} \subset A$.
\end{thm}

Recently, Crann et al. have verified a conjecture in \cite{BG} (see Conjecture $3.20$) using \Cref{kribs} as the main technical tool. Indeed, \Cref{kribs} has the following marginal generalization.

\begin{ppsn}\label{bhat}
If $B\subset A$ is a connected inclusion of finite-dimensional $C^*$-algebra with $B$ simple, then $A/B$ has a unitary orthonormal basis with respect to the Markov trace.
\end{ppsn}
\begin{prf}
Since $B$ is simple, $B\subset A$ is isomorphic to $B\subset B\otimes (B^{\prime}\cap A)$. The conclusion now follows from \Cref{kribs}.\qed
\end{prf}

\noindent\textbf{Proof of \Cref{to prove}\,:} Suppose that $N\subset M$ is a subfactor with $N^{\prime}\cap M$ is simple, that is, $N^{\prime}\cap M\cong M_n(\mathbb{C})$ for some $n\in\bbn$. 
\smallskip

\textbf{Step $(i)$:} We first show that, $M_1 \subset M_2$ has a unitary orthonormal basis. Note that $N'\cap M \ni x \mapsto JxJ \in M'\cap M_1$ is an anti-isomorphism; so, $M'\cap M_1$ is also simple. Further,  since $N \subset M$ has depth $2$, it follows that  $M \subset M_1$ is also of depth 2. Thus, by \Cref{wellknown} we know that $tr_{M_2}|_{M^{\prime}\cap M_2}$ is the Markov trace for $M'\cap M_1 \subset M'\cap M_2$. By \Cref{bhat}, $M'\cap M_2$ admits a unitary orthonormal basis , say  $\{u_i : 1 \leq i \leq d\}$ over $M^{\prime}\cap M_1$, w.r.t. this trace. Moreover, we know that
\[
\begin{array}{ccc}
M_1 & \subset & M_2\\
\cup & & \cup\\
M'\cap M_1 & \subset & M'\cap M_2\end{array}
\]
is a non-degenerate commuting square. Hence, $\{u_i\}$ is also an orthonormal unitary basis for $M_2/M_1$ (see \Cref{combasis}).  
\smallskip

\textbf{Step $(ii)$:} We prove that $N \subset M$ has a unitary orthonormal basis. Indeed, fix a $2$-step downward basic construction $N_{-2}\subset N_{-1}\subset N\subset M$ and recall that (see for instance, \cite[Prop. 3.12]{BG2}) that $N_{-2}\subset N_{-1}$ has depth $2$. Also, $N_{-2}\subset N_{-1}$ has simple first relative commutant because $N_{-2}^\prime\cap N_{-1}\cong N^\prime\cap M$ (via the shift operator - see the proof of Cor. $3.12$, \cite{BG2}). Hence, using Step $(i)$ for the subfactor $N_{-2} \subset N_{-1}$, we conclude that  $N \subset M$ has a unitary orthonormal basis.
\smallskip

These two steps together completes the proof. \qed

%%%%%%%%%%%%%%%%%%%%%%%%%%%%%%%%%%%%%%%
%%%%%%%%%%%%%%%%%%%%%%%%%%%%%%%%%%%%%%%
%%%%%%%%%%%%%%%%%%%%%%%%%%%%%%%%%%%%%%%

\newsection{On Popa's relative dimension of projections}\label{Sec 3}

In \cite{popa}, Popa introduced the set $\Lambda(M,N)$ to measure the dimension of projections $p\in M$ relative to the subfactor $N$. Computation of $\Lambda(M,N)$ in full generality is extremely hard and remains an open problem (see \cite{popa, pop2}). In \cite{popa}, this has been viewed as an orthogonalization problem (see also \cite{pop2}). More recently, in \cite{BDLR}, a notion of angle between interemdiate subfactors of a subfactor with finite Jones index have been discovered. In \cite{BG1}, the possible value of angle has been shown to be intimately related with Popa's question of finding the possible elements in $\Lambda$. 
\smallskip

\noindent\textbf{Notation $(3)$:} For a von Neumann algebra $\mathcal{M}$, we denote by $\mathcal{P(M)}$ the set of all projections in $\mathcal{M}$.
\begin{dfn}[\cite{popa}]
For a subfactor $N\subset M$, define $\Lambda(M,N)=\{\alpha\geq 0\,:\,E^M_N(p)=\alpha,\,p\in\mathcal{P}(M)\}$.
\end{dfn}

Observe that $0,1\in\Lambda(M,N)$, and since $E$ is a contraction, $\Lambda(M,N)\subseteq[0,1]$. It is known that $[M:N]^{-1}=\inf\,\Lambda(M,N)$ \cite{popa}. The following results will be repeatedly used possibly without any further mention.

\begin{ppsn}[Proposition $1.7$, \cite{popa}]
If $(\mathcal{N\subset Q,P\subset M})$ is a commuting square of von Neumann algebras, then $\Lambda(\mathcal{Q,N})\subseteq\Lambda(\mathcal{M,P})$.
\end{ppsn}

\begin{ppsn}[Proposition $4.1$, \cite{popa}]\label{e}
Let $N\subset M$ be a finite index subfactor with basic construction $M\subset M_1$. For any $\beta>0$, let $M^\beta$ be the $\beta$-amplification of $M$. Then, we have the following\,:
\begin{enumerate}[$(i)$]
\item If $\alpha\in\Lambda(M,N)$, then $1-\alpha\in\Lambda(M,N)$.
\item If $\alpha\in\Lambda(M,N),\,\alpha\neq 0,1$, then $\lambda/\alpha\in\Lambda(M^\alpha_1,M^\alpha)$ and $\lambda/(1-\alpha)\in\Lambda(M_1^{1-\alpha},M^{1-\alpha})$, where $\lambda=[M:N]^{-1}$.
\end{enumerate}
\end{ppsn}

\begin{lmma}\label{e0}
For any unitary $u\in\mathcal{U}(M)$, we have $\Lambda(M,N)=\Lambda(M,uNu^*)$.
\end{lmma}
\begin{prf}
For $\alpha\in\Lambda(M,N)$, there is a projection $p\in M$ such that $E(p)=\alpha$. Consider $\mbox{Ad}_u(p)\in M$. Then,
\[
E_{uNu^*}(\mbox{Ad}_u(p))=\mbox{Ad}_u\circ E^M_N(p)=\alpha\,.
\]
Hence, $\Lambda(M,N)\subseteq\Lambda(M,uNu^*)$. Replacing $u$ by $u^*$ in the above, similar analysis gives the reverse inclusion.\qed
\end{prf}

The set $\Lambda(M,N)$ is completely known when $[M:N]\leq 4$ and $N$ is locally trivial; and only some parts of $\Lambda(M,N)$ is known when $[M:N]>4$. Let $0<t<\frac{1}{2}$ be the real number given by the quadratic equation $t(1-t)=[M:N]^{-1}$. For $[M:N]\leq 4$, both the sets $\Lambda(M,N)\cap[0,t]$ and $\Lambda(M,N)\cap[1-t,1]$ are computed in \cite{popa}. Note that $\Lambda$ depends only on the index $[M:N]$ in all these cases. However, for index greater than $4$, the set $\Lambda(M,N)\cap(t,1-t)$ remains mysterious. In \cite{popa}, one of the major questions asked by Popa is whether the set $\Lambda(M,N)\cap(t,1-t)$ is nonempty, at least for irreducible subfactors, where $t(1-t)=[M:N]^{-1}$ and $t<1/2$ (see $5.4.3$ in \cite{popa}). Goal of this section is to show that at least for three different classes of subfactors, this question of Popa has positive answer.

Before we begin, let us get familiar with a construction of Popa in Proposition $5.5$, \cite{popa}. If $\Lambda(M,N)\cap(t,1-t)$ has one element $\alpha$ such that $t<\alpha\leq 1/2$, then in certain situations, it is possible to construct a sequence of elements in the set starting with $\alpha$. If $N\subset M$ is such that $[M:N]>4$ and the inclusion $N\subset M$ splits $R$ (i.e., there exists $R\subset N$ with $N=R\vee(R^\prime\cap N)$ and $M=R\vee(R^\prime\cap M)$), then for any $t<\alpha< 1/2$, there exists a sequence $\{\alpha_m\}_{m\geq 0}$ defined recursively by $\alpha_0=\alpha$ and $\alpha_{m+1}=\frac{[M:N]^{-1}}{1-\alpha_m}$ such that $\{\alpha_{2m}\}_{m\geq 0}\subseteq\Lambda(N\subset M)$ and $\{\alpha_{2m+1}\}_{m\geq 0}\subseteq\Lambda(M\subset M_1)$. Moreover, $\alpha_m$ decreases to $t$. Thus, if $N\subset M$ satisfies the above hypothesis, and if $\Lambda(M,N)\cap(t,1-t)$ is nonempty, then the cardinality $\texttt{\#}\,\Lambda(M,N)\cap(t,1-t)$ is at least $\aleph$.\smallskip

Recall that if $M/N$ has a unitary orthonormal basis, then $[M:N]$ must be an integer.

\begin{thm}\label{popa}
If $M/N$ has a unitary orthonormal basis then $\{0, 1/n, 2/n, \cdots, (n-1)/n,1\}\subseteq \Lambda(M_1,M)$, where $[M:N]=n$.
\end{thm}
\begin{prf}
Suppose $\{u_i:1\leq i\leq n\}$ is a unitary orthonormal basis. Therefore, $E_N(u^*_iu_j)=\delta_{ij}$.  For $1\leq m\leq n$ we see that
$\big(\sum_{i=1}^m u_ie_1u^*_i\big)\big(\sum_{j=1}^m u_je_1u^*_j\big)=\sum_{i,j=1}^m u_iE_N(u^*_iu_j)e_1u^*_j=\sum_{i=1}^m u_ie_1u^*_i$. Thus, for each $1\leq m\leq n$, the self-adjoint operator
$\sum_{i=1}^m u_ie_1u^*_i$ is a projection in $M_1.$ Also, $E^{M_1}_M(\sum_{i=1}^m u_ie_1u^*_i)=m/n$, which completes the proof.\qed
\end{prf}

\begin{crlre}\label{depth2}
Suppose that $N\subset M$ is a depth $2$ subfactor with $N^{\prime}\cap M$ simple. Then, $[M:N]$ is an integer, say $n$, and $\{0,1/n, 2/n,\cdots, (n-1)/n,1\}\subseteq \Lambda(M,N)$.
\end{crlre}
\begin{prf}
A subfactor $N\subset M$ is of depth 2 if and only if the subfactor $N_{-1}\subset N$ is of depth 2, where $N_{-1}\subset N\subset M$ is a (semi-canonical) tower of downward basic construction. As $N^{\prime}_{-1}\cap N$  is anti-isomorphic to $N^{\prime}\cap M$,  we find that $N^{\prime}_{-1}\cap N$ is simple. The rest follows from \Cref{to prove} and \Cref{popa}. \qed
\end{prf}

\begin{crlre}
Let $N\subset M$ be a subfactor with $[M:N]=n\in \mathbb{N}$. If $M/N$ has a two-sided orthonormal basis, then $\{0,1/n, 2/n, \cdots, (n-1)/n,1\}\subseteq \Lambda(M_1,N).$
\end{crlre}
\begin{prf}
Suppose that $\{\lambda_1,\cdots, \lambda_n\}$ is an orthonormal basis for $M/N$, and thus $\sum_{i=1}^n \lambda_ie_1\lambda^*_i=1$. As in the proof of \Cref{popa}, we see that for each $1\leq m\leq n$, the operator $\sum_{i=1}^m \lambda_ie_1\lambda^*_i$ is a projection in $M_1$. Now,
\[
E^{M_1}_N\big(\sum_{i=1}^m \lambda_ie_1\lambda^*_i\big)=E^{M}_N\circ E^{M_1}_M\big(\sum_{i=1}^m \lambda_ie_1\lambda^*_i\big)=\frac{1}{n}\sum_{i=1}^mE_N(\lambda_i\lambda^*_i)=m/n,
\]
where the last equality follows from the fact that $\{\lambda^*_i:1\leq i\leq n\}$ is also an orthonormal basis, if $\{\lambda_i:1\leq i\leq n\}$ is so.\qed
\end{prf}

\begin{thm}\label{depth2popa}
For a depth $2$ subfactor $N\subset M$ such that $[M:N]>4$ and $N^\prime\cap M$ is simple, we have $\Lambda(M,N)\cap(t,1-t)\neq\emptyset$, where $0<t<\frac{1}{2}$ is given by the quadratic equation $t(1-t)=[M:N]^{-1}$.
\end{thm}
\begin{prf}
For a depth $2$ subfactor $N\subset M$ with $N^\prime\cap M$ is simple, we have $[M:N]=n$, an integer, and $\{0,1/n,\ldots,(n-1)/n,1\}\subset\Lambda(M,N)$ by \Cref{depth2}. Consider $k=\lfloor\frac{n+1}{2}\rfloor$. Then, $k/n\in\Lambda(M,N)$. If $0<t< 1/2$ and $t(1-t)=1/n$ with $n>4$, then $t=\frac{1}{2}\left(1-(1-\frac{4}{n})^{\frac{1}{2}}\right)$. Now, in both the situations when $n$ is even and odd, we see that $k/n>t$. This is because, on contrary, if we assume that $k/n\leq t$, then we get that $\lfloor\frac{n+1}{n}\rfloor=1\leq\left(1-(1-\frac{4}{n})^{\frac{1}{2}}\right)$, a contradiction to $n>4$. Similarly, we can conclude that $k/n<1-t$, that is $k/n\in(t,1-t)$. Hence, $k/n\in\Lambda(M,N)\cap(t,1-t)$.\qed
\end{prf}

\begin{ppsn}\label{spin}
Let $u$ be an $n\times n$ complex Hadamard matrix where $n>4$, and $R_u\subset R$ be the associated spin model subfactor of index $\,n$. Then, $\big\{\frac{k}{n}\,:\,0\leq k\leq n\big\}\subseteq\Lambda(R,R_u)$. Moreover, $\Lambda(R,R_u)\cap(t,1-t)\neq\emptyset$, where $0<t<\frac{1}{2}$ is given by the quadratic equation $t(1-t)=\frac{1}{n}$.
\end{ppsn}
\begin{prf}
Recall from \Cref{preliminaries} the following non-degenerate commuting square
\[
\begin{matrix}
\mbox{Ad}_u(\Delta_n) &\subset & M_n\\
\cup & &\cup\\
\bbc &\subset & \Delta_n
\end{matrix}
\]
where $u$ is an $n\times n$ complex Hadamard matrix, and iterating Jones' basic construction vertically, the spin model subfactor $R_u\subset R$ is obtained with $[R:R_u]=n$. Therefore, we have the following commuting square
\[
\begin{matrix}
R_u & \subset & R\\
\cup & & \cup\\
\bbc & \subset & \Delta_n\\
\end{matrix}
\]
and hence, $\Lambda(\Delta_n,\bbc)\subseteq\Lambda(R,R_u)$. The conditional expectation $E:\Delta_n\to\bbc$ is given by the unique normalized trace $tr$ on $M_n$. Since there are projections of rank $k$ in $\Delta_n$, where $1\leq k\leq n$, it is now immediate that $\{0,1/n,\ldots,(n-1)/n,1\}\subseteq\Lambda(\Delta_n,\bbc)$. Hence, we have $\{0,1/n,\ldots,(n-1)/n,1\}\subseteq\Lambda(R,R_u)$. Now, the similar trick used in the proof of \Cref{depth2popa} gives the result.\qed
\end{prf}

\begin{thm}\label{aa}
Let $u$ be an $n\times n$ complex Hadamard matrix and $n>4$. For the spin model subfactor $R_u\subset R$, we have
\[
\left\{\frac{1}{k}: 1\leq k\leq n\right\}\cup\left\{\frac{k-1}{k}:1\leq k\leq n\right\}\subseteq\Lambda(R,R_u).
\]
Moreover, if $\,2\leq k\leq n-2$, then each $\frac{1}{k}$ and $\frac{k-1}{k}$ lie in $\Lambda(R,R_u)\cap(t,1-t)$, where $t(1-t)=\frac{1}{n}$ and $t>1/2$.
\end{thm}
\begin{prf}
Suppose that $R_{u,-1}\subset R_u\subset R$ is an instance of downward basic construction. By Proposition $4$ in \cite{KSV}, we see that $R_u/R_{u,-1}$ has a unitary orthonormal basis. Thus, we have $\frac{m}{n}\in\Lambda(R_{u,-1}\subset R_u)$ for $0\leq m\leq n$. Consider $\alpha=\frac{k}{n}$, where $0<k<n$. By Proposition $4.1$ in \cite{popa}, we have $\frac{1}{k}\in\Lambda(R^\alpha,R_u^\alpha)=\Lambda(R,R_u)$. Now, if $\alpha\in\Lambda(R,R_u)$, then $1-\alpha\in\Lambda(R,R_u)$, and moreover, both $1$ and $\frac{1}{n}$ belong to $\Lambda(R,R_u)$ by \Cref{spin}. Hence, the first part follows. To see the second part, first observe that $t<\frac{1}{k}<1-t$ if and only if $\frac{k^2}{k-1}<n$. Moreover, $t<\frac{1}{k}<1-t$ implies that $t<\frac{k-1}{k}<1-t$. Now, the inequality $\frac{k^2}{k-1}<n$ holds for any $1<k\leq n-3$ because
\[
k^2\leq (k-1)(k+2)\Rightarrow \frac{k^2}{k-1}\leq k+2\Rightarrow \frac{k^2}{k-1}\leq n-1<n
\]
when $k>1$. Also, the inequality holds for $k=n-2$ as $n>4$, and it never holds for $k=n,n-1$. This completes the second part.\qed
\end{prf}

\begin{ppsn}\label{vertex}
Let $v$ be an $n^2\times n^2$ bi-unitary matrix, where $n>2$, and $R_v\subset R$ be the corresponding vertex model subfactor of index $\,n^2$. Then, $\big\{\frac{k}{n}\,:\,0\leq k\leq n\big\}\subseteq\Lambda(R,R_v)$. Moreover, $\Lambda(R,R_v)\cap(t,1-t)\neq\emptyset$, where $0<t<\frac{1}{2}$ is given by the quadratic equation $t(1-t)=\frac{1}{n^2}$.
\end{ppsn}
\begin{prf}
Recall from \Cref{preliminaries} the following non-degenerate commuting square
\[
\begin{matrix}
\mbox{Ad}_v(\bbc\otimes M_n) &\subset & M_n\otimes M_n\\
\cup & &\cup\\
\bbc &\subset & M_n\otimes\bbc
\end{matrix}
\]
where $v$ is an $n^2\times n^2$ bi-unitary matrix, and iterating Jones' basic construction vertically, the vertex model $R_v\subset R$ is obtained with $[R:R_v]=n^2$. Therefore, we have the following commuting square
\[
\begin{matrix}
R_v & \subset & R\\
\cup & & \cup\\
\bbc & \subset & M_n\otimes\bbc\\
\end{matrix}
\]
and hence, $\Lambda(M_n,\bbc)\subseteq\Lambda(R,R_v)$. The conditional expectation $E:M_n\to\bbc$ is given by the unique normalized trace $tr$ on $M_n$. Since there are projections of rank $k$ in $M_n$, where $1\leq k\leq n$, it is now immediate that $\{0,1/n,\ldots,(n-1)/n,1\}\subseteq\Lambda(M_n,\bbc)$. The rest of the proof is similar to that of \Cref{spin} and we skip the details.\qed
\end{prf}

\begin{thm}\label{bb}
Let $v$ be an $n^2\times n^2$ bi-unitary matrix, where $n>2$. For the vertex model subfactor $R_v\subset R$, we have
\[
\left\{\frac{1}{nk}:1\leq k\leq n\right\}\cup\left\{\frac{nk-1}{k}:1\leq k\leq n\right\}\subset\Lambda(R,R_v).
\]
Moreover, for $1\leq k\leq n-1$ each $\frac{1}{nk}$ and $\frac{nk-1}{k}$ lie in $\Lambda(R,R_u)\cap(t,1-t)$, where $t(1-t)=\frac{1}{n^2}$ and $t>1/2$.
\end{thm}
\begin{prf}
Suppose that $R_{v,-1}\subset R_v\subset R$ is an instance of downward basic construction. Again by Proposition $4$ in \cite{KSV}, we see that $R_v/R_{v,-1}$ has a unitary orthonormal basis. Thus, by \Cref{vertex} we have $\frac{m}{n}\in\Lambda(R,R_v)$ for $0\leq m\leq n$. Consider $\alpha=\frac{k}{n}$, where $1<k<n$. By Proposition $4.1$ in \cite{popa}, we have $\frac{\lambda(R,R_v)}{\alpha}=\frac{1}{nk}\in\Lambda(R^\alpha,R_v^\alpha)=\Lambda(R,R_v)$. Now, if $\alpha\in\Lambda(R,R_v)$, then $1-\alpha\in\Lambda(R,R_v)$, and moreover both $1$ and $\frac{1}{n^2}$ belong to $\Lambda(R,R_v)$. Hence, the first part follows. For the second part, first observe that $t<\frac{1}{nk}<1-t$ if and only if $\frac{k^2+1}{k}<n$, and moreover $t<\frac{1}{nk}<1-t$ implies that $t<\frac{nk-1}{k}<1-t$. Now, the inequality $\frac{k^2+1}{k}<n$ holds for any $1<k<n$ because
\[
k\leq n-1\Rightarrow k^2\leq(n-1)k\Rightarrow k^2+1\leq nk+(1-k)\Rightarrow k^2+1<nk
\]
when $k>1$. Since $n>2$, the inequality holds for $k=1$, and it never holds for $k=n$. This completes the second part.\qed
\end{prf}

\begin{rmrk}\rm
We end this section with the remark that for the spin model and vertex model subfactors, using Popa's construction discussed at the beginning of the section (or see Proposition $5.5$, \cite{popa}), we get countable infinitely many elements in the set $\Lambda\cap(t,1-t)$. This is because if $N\subset M$ is either a spin model or a vertex model subfactor, then the inclusion $N\subset M$ splits $R$ (i.e., there exists $R\subset N$ with $N=R\vee(R^\prime\cap N)$ and $M=R\vee(R^\prime\cap M)$), since $N\subset M$ is obtained as iterated basic construction of non-degenerate commuting squares of finite-dimensional algebras. Therefore, for each element $\alpha=\frac{1}{k}$ in \Cref{aa}, and $\beta=\frac{1}{nk}$ in \Cref{bb}, for the cases of spin model subfactor $R_u\subset R$ and vertex model subfactor $R_v\subset R$ respectively, we get sequences $\{\alpha_{2m}\}_{m\geq 0}$ in $\Lambda(R_u\subset R)\cap(t,1-t)$ and $\{\beta_{2m}\}_{m\geq 0}$ in $\Lambda(R_v\subset R)\cap(t,1-t)$.
\end{rmrk}

%%%%%%%%%%%%%%%%%%%%%%%%%%%%%%%%%%%%%
%%%%%%%%%%%%%%%%%%%%%%%%%%%%%%%%%%%%%
%%%%%%%%%%%%%%%%%%%%%%%%%%%%%%%%%%%%%

\newsection{A construction of infinitely many new elements in $\Lambda$ for the cases of spin and vertex model subfactors}\label{Sec 4}

In \Cref{Sec 3}, we have shown that if $N\subset M$ is a spin model or a vertex model subfactor, then $\Lambda(N\subset M)\cap(t,1-t)\neq\emptyset$, where $t\leq\frac{1}{2}$ is given by $t(1-t)=[M:N]^{-1}$. In fact, we have produced finitely many elements in the intersection, and then applying Popa's construction one can produce a sequence of elements in $\Lambda(N\subset M)\cap(t,1-t)$. Goal of this section is to construct infinitely many elements (distinct from those arising through Popa's construction) in $\Lambda(R_v\subset R)\cap(t,1-t)$, where $R_v\subset R$ is either an even index spin model subfactor or a vertex model subfactor of index $(2n)^2,\,n\geq 3$. Then, applying Popa's construction to each element in our sequence, we can produce a lot of new elements in $\Lambda(N\subset M)\cap(t,1-t)$ for these two classes of subfactors. We begin by fixing some notations that are crucial in this section.
\smallskip

\noindent\textbf{Notation $(4)$:} Let $k\in\mathbb{N}\cup\{0\}$ and $u$ be a comple Hadamard matrix of order $n$. Consider the following unitary matrices\,:
\[
D_u:=\sqrt{n} \sum_{i=1}^n \sum_{j=1}^n \overline{u}_{ij} (E_{ii} \otimes E_{jj})
\] 
and
\begin{IEEEeqnarray}{lCl}\label{Not5}
u_{2k+1} &:=& (I_n \otimes u_{2k})(D_u \otimes I_n^{(k)})\nonumber\\
u_{2k} &:=& u_{2k-1}(u \otimes I_n^{(k)})
\end{IEEEeqnarray}
defined recursively with the convention $u_0=u$. Note that these unitary matrices appear in the grid of finite-dimensional algebras for the spin model subfactor $R_u\subset R$ (see \Cref{preliminaries}).
\smallskip

\begin{figure} 
\centering
\begin{tikzpicture}[font=\large,thick]

\node[draw,
    align=center,
    minimum width=2.5cm,
    minimum height=1cm] (block1) { Take any spin model ~\\
     subfactor $\,R_u\subset R$ of index $4$ };

\node[draw,
    align=center,
    below=of block1,
    minimum width=3.5cm,
    minimum height=1cm] (block2) {~Construct a sequence $\{\alpha_m\}\subseteq\Lambda(R_u\subset R)\,$\\
     that comes from the f.d. grid for $R_u\subset R\,$~};

\node[draw,
    align=center,
    right=of block2,
    minimum width=3.5cm,
    minimum height=1cm] (block3) {~Take any even index, say $n$, spin ~\\
     model subfactor $R_v\subset R$~};
  
\node[draw,
    align=center,
    below=of block3,
    minimum width=3.5cm,
    minimum height=1cm] (block4) {~Construct $n-3$ injective\\ functions $f_j$ such that\\ $\{f_j(\alpha_m):j=1,\ldots,n-3;\,m\in\bbn\}\subseteq$\\ $\Lambda(R_v\subset R)\cap (t,1-t)$~};

\node[draw,
    rectangle,
    left=of block4,
    minimum width=3.5cm,
    minimum height=1cm] (block5) {~ Apply similar scheme for
     vertex model ~};    

\draw[-latex] (block1) edge (block2)
    (block2) edge (block3)
    (block3) edge (block4)
    (block4) edge (block5);
    
\end{tikzpicture}
\bigskip

\caption{Scheme for constructing infinitely many elements in $\Lambda(R_v\subset R)\cap(t,1-t)$}\label{scheme}
\end{figure}
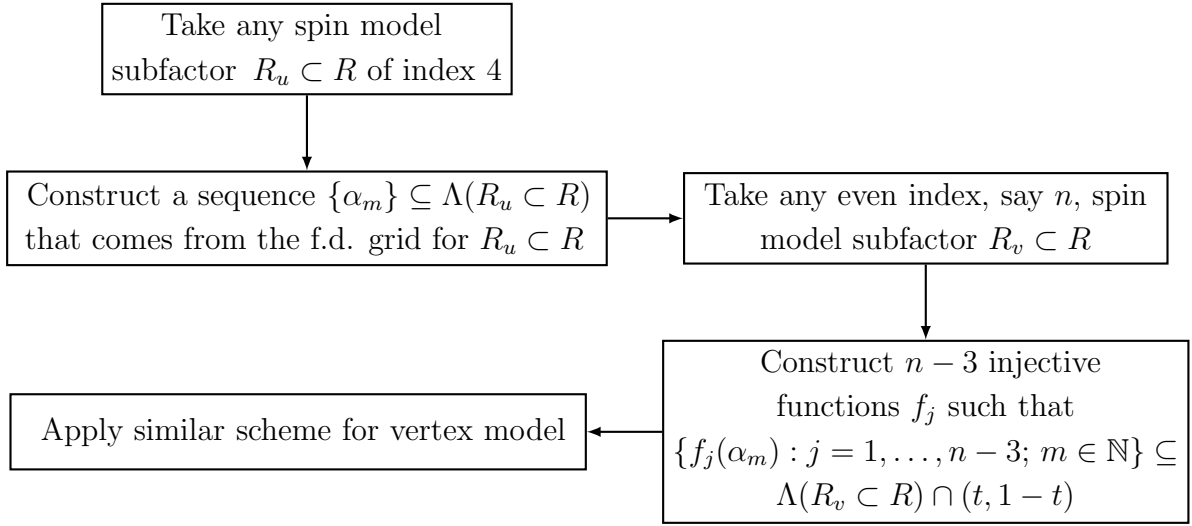

\noindent\textbf{Our scheme:} We first consider any spin model subfactor $R_u\subset R$ of index $4$, that is, $u$ is a $4\times 4$ complex Hadamard matrix. Due to \cite{popa}, we know all the elements in $\Lambda(R_u\subset R)$. We identify a sequence $\{\alpha_m\}\subseteq\Lambda(R_u\subset R)$ such that each $\alpha_m$ appears from the grid $\mbox{Ad}_{u_{2k-1}}(M_4^{(k)})\subset \Delta_4\otimes M_4^{(k)},\,k\in\bbn$, of finite-dimensional algebras for the subfactor $R_u\subset R$ (see \Cref{Not5} in this regard). That is, for each $m\in\bbn$,
\[
\alpha_m\in\Lambda\big(M_4^{(k)})\subset \Delta_4\otimes M_4^{(k)}\big)\subseteq\Lambda(R_u\subset R)
\]
for some $k\in\bbn$. Now, consider any even index spin model subfactor $R_v\subset R$, where $v$ is an $n\times n$ complex Hadamard matrix such that $n=[R:R_v]\geq 6$ is even integer. We construct $n-3$ injective functions $f_j$ and show that $\{f_j(\alpha_m):j=1,\ldots,n-3;\,m\in\bbn\}\subseteq\Lambda(R_v\subset R)\cap(t,1-t)$, where $t<\frac{1}{2}$ and $t(1-t)=1/n$. It turns out that this construction can be carried to the case of vertex model subfactors also. This scheme is depicted in \Cref{scheme}.
\medskip

Let us now proceed to the construction. Let $R_u\subset R$ be any spin model subfactor of index $4$, and consider any spin model subfactor $R_v\subset R$ of index $n$ such that $n\in 2\bbn$ and $n\geq 6$. Recall that the alternating stages in the finite-dimensional grid for the subfactor $R_v\subset R$ is of the form $\mbox{Ad}_{v_{2k-1}}(\bbc\otimes M_{n^k})\subset \Delta_n\otimes M_{n^k}$, where $k\in\bbn$ and $v_{2k-1}\in\Delta_n\otimes M_{n^k}$ are unitary matrices (see \Cref{Not5} in this regard). For each $k\in\bbn$, the conditional expectation $E_k:\Delta_n\otimes M_{n^k}\to\mbox{Ad}_{v_{2k-1}}(\bbc\otimes M_{n^k})$ is given by $\mbox{Ad}_{v_{2k-1}}(\mbox{tr}\otimes\mbox{id})\mbox{Ad}_{v_{2k-1}^*}$, where $\mbox{tr}:\Delta_n\to\bbc$ is the restriction of the unique normalized trace on $M_n$. Since
\[
\Lambda(\mbox{Ad}_{v_{2k-1}}(\bbc\otimes M_{n^k})\subset \Delta_n\otimes M_{n^k})=\Lambda(\bbc\otimes M_{n^k}\subset \Delta_n\otimes M_{n^k})
\]
due to \Cref{e0}, if we are interested only in $\Lambda$, then we can safely discard the unitary matrices $v_{2k-1}$ at each stage and $E_k:\Delta_n\otimes M_{n^k}\to M_{n^k}$ is given by $(x_1,\ldots,x_n)\mapsto\frac{1}{n}(x_1+\cdots+x_n)$, where $x_j\in M_{n^k}$. Note that
\[
\Lambda(M_{n^k}\subset \Delta_n\otimes M_{n^k})\subseteq\Lambda(M_{n^{k+1}}\subset \Delta_n\otimes M_{n^{k+1}})\,,
\]
and consequently
\[ 
\cup_k\Lambda(M_{n^k}\subset \Delta_n\otimes M_{n^k})\subseteq\Lambda(R_u\subset R)
\]
due to the fact that
\[
\begin{matrix}
\mbox{Ad}_{v_{2k+1}}(\bbc\otimes M_{n^{k+1}}) &\subset & \Delta_n\otimes M_{n^{k+1}}\\
\cup & &\cup\\
\bbc\otimes \mbox{Ad}_{v_{2k-1}}(\bbc\otimes M_{n^k}) &\subset & \bbc\otimes\Delta_n\otimes M_{n^k}\\ 
\end{matrix}
\]
are commuting squares by the construction of $R_v\subset R$. The set $\cup_k\Lambda(M_{n^k}\subset \Delta_n\otimes M_{n^k})$ consists of precisely those $\alpha\in\Lambda(R_u\subset R)$ that arise from the finite-dimensional grid for $R_v\subset R$.

Any projection $p\in\Delta_n\otimes M_{n^k}$ is given by $n$-tuple $(p_1,\ldots,p_n)$ of projections $p_j\in M_{n^k},\,j=1,\ldots,n$. Therefore, by definition of $\Lambda$, we get that $\alpha\in\Lambda(M_{n^k}\subset \Delta_n\otimes M_{n^k})$ if and only if $\frac{1}{n}(\sum_{j=1}^np_j)=\alpha I_{n^k}$. Note that some of the $p_j$ can be zero also. Thus, solution to the equation $E_k(p)=\alpha I_{n^k}$ for $p\in\mathcal{P}(\Delta_n\otimes M_{n^k})$ reduces to the existence of at most $n$-projections $p_j,\,j=1,\ldots,n,$ in $M_{n^k}$ such that
\[
p_1+\cdots+p_n=n\alpha I_{n^k}\,.
\]
This is the same as asking whether the following unital $\star$-algebra
\[
\mathscr{A}_{n,\beta}:=\bbc\,\Big\langle p_1,\ldots,p_n\,:\,p_j^*=p_j=p_j^2,\,j=1,\ldots,n,\,\sum_{j=1}^np_j=\beta e\Big\rangle\,,
\]
where $e$ denotes the unit element, has at least one $\star$-representation on the finite-dimensional Hilbert spaces $\mathcal{H}_n(k):=\bbc^{n^k}$ for $k\in\bbn$. At this point, it is worthwhile to mention that finding values of the parameter $\tau\in\mathbb{R}$ such that the $\star$-algebra $\mathcal{T}\mathcal{L}_{\infty,\tau}=\bbc\big\langle p_1,\ldots,p_n,\ldots\,:\,p_k^2=p_k=p_k^*;\,p_kp_j=p_jp_k,\,|k-j|\geq 2;\,p_kp_{k\pm 1}p_k=\tau p_k\big\rangle$ has at least one representation is in similar theme and goes back to the celebrated work of Jones \cite{J}.

In this regard, we mention that motivated by \cite{K, Fu}, numerous authors have investigated representation of the $\star$-algebra $\mathscr{A}_{n,\beta}$ using the Coxeter functors (see e.g., \cite{F, RS, RS1, KRS1, KRS, S}, and references therein). However, there was no restriction on the dimensions of the Hilbert spaces on which representations exist, and moreover, infinite-dimensional Hilbert spaces were also allowed. This is where our situation differs. We need to study whether representation of $\mathscr{A}_{n,\beta}$ exists on $n^k$-dimensional ($k\in\bbn$) Hilbert spaces. Thus, our situation is drastically much more complicated. Note that since trace of a finite-dimensional projection is rational, $\beta$ can take only rational values in any representation of $\mathscr{A}_{n,\beta}$ on $\mathcal{H}_n(k)$. It is worthwhile to investigate all possible $\star$-representations of $\mathscr{A}_{n,\beta}$ on $\mathcal{H}_n(k),\,k\in\bbn$.

Consider the $n=4$ case. A projection in $\Delta_4\otimes M_{4^k}$ is a tuple $(p_1,p_2,p_3,p_4)$ such that each $p_j\in M_{4^k}$. By definition of $\Lambda$, we get that $\alpha\in\Lambda(M_{4^k}\subset \Delta_4\otimes M_{4^k})$ if and only if $\frac{1}{4}(\sum_{j=1}^4p_j)=\alpha I_{4^k}$. Note that some of the $p_j$'s can be zero also. Thus, we are interested in the existence of at most four projections $p_j,\,j=1,2,3,4,$ in $M_{4^k}$ such that
\begin{IEEEeqnarray}{lCl}\label{ux}
p_1+p_2+p_3+p_4=4\alpha I_{4^k}\,.
\end{IEEEeqnarray}
In other words, we concentrate on $\star$-representations of the unital $\star$-algebra $\mathscr{A}_{4,\beta}$ on the finite dimensional Hilbert spaces $\mathcal{H}_4(k):=\bbc^{4^k}$ for $k\in\bbn$. Consider the set $\Sigma_4=\{2\}\cup\wedge_4\cup(4-\wedge_4)$ (see \cite{KRS}), where
\[
\wedge_4:=\{0,\Phi_4(0),\Phi_4(\Phi_4(0)),\ldots\}\cup\{1,\Phi_4(1),\Phi_4(\Phi_4(1)),\ldots\}
\]
and
\begin{IEEEeqnarray}{lCl}\label{pq}
\Phi_4(x):=1+\frac{1}{4-x-1}=1+\frac{1}{3-x}\,.
\end{IEEEeqnarray}
Note that $\Sigma_4$ contains the set of $\beta$ such that $\mathscr{A}_{4,\beta}$ has $\star$-representation on a finite dimensional Hilbert space. Now, recall from \cite{S}, the following definitions of $\Sigma_n\mbox{ and }\wedge_n$ for general $n\geq 5\,:$
\[
\wedge_n:=\wedge_n^1\cup\wedge_n^2,
\]
\[
\wedge_n^1:=\Big\{0\,,\,1+\frac{1}{n-1}\,,\,1+\frac{1}{n-2-\frac{1}{n-1}}\,,\,1+\frac{1}{n-2-\frac{1}{n-2-\frac{1}{n-1}}}\,,\,\ldots\Big\},
\]
\[
\wedge_n^2:=\Big\{1\,,\,1+\frac{1}{n-2}\,,\,1+\frac{1}{n-2-\frac{1}{n-2}}\,,\,1+\frac{1}{n-2-\frac{1}{n-2-\frac{1}{n-2}}}\,,\,\ldots\Big\},
\]
\[
\Sigma_n:=\wedge_n\cup\left[\frac{n-\sqrt{n^2-4n}}{2},\frac{n+\sqrt{n^2-4n}}{2}\right]\cup(n-\wedge_n).
\]
Notice that when $n\geq 5$, the set $\Sigma_n$ contains an interval, and $\wedge_n$'s are discrete. Also, for $t<\frac{1}{2}$ satisfying the quadratic equation $t(1-t)=\frac{1}{n}$, the interval in $\Sigma_n$ is precisely $[nt,n-nt]$. Our purpose of recalling the above sets is the following result.

\begin{thm}[\cite{KRS}]\label{known}
A real matrix $\alpha I_n$ can be written as sum of $r$ projections if and only if $\alpha\in\Sigma_r$ and $n\alpha\in\bbn$.
\end{thm}
We apply this theorem for $r=4$ and $n=4^k$ for $k\in\bbn$. Thus, we have a solution to the equation $\sum_{j=1}^4p_j=\beta I_{4^k}$ in \Cref{ux} if and only if $\beta\in\Sigma_4$ and $4^k\beta\in\bbn$. Now, consider the subset $\{1,\Phi_4(1),\Phi_4(\Phi_4(1)),\ldots\}$ of $\wedge_4\subseteq\Sigma_4$.

\begin{lmma}\label{s}
For any $k\in\bbn$, we have $\Phi_4^{(k)}(1)=\frac{3+2(k-1)}{2+(k-1)}$, where $\Phi_4^{(k)}(1)$ is defined recursively as $\Phi_4\big(\Phi_4^{(k-1)}(1)\big)$.
\end{lmma}
\begin{prf}
This follows easily by induction on $k$, using the formula for $\Phi_4(x)$ described in the \Cref{pq}.\qed
\end{prf}

It is now obvious that $4^r\Phi_4^{(k)}(1)\in\bbn$ for some $r,k\in\bbn$ if and only if $(k+1)$ divides $4^r$, that is, $k$ must be of the form $2^m-1$ for $m\in\bbn$.

\begin{ppsn}\label{qp}
Let $\beta_m=\Phi_4^{(2^m-1)}(1)\in\wedge_4$ for $m\in\bbn$, and consider a spin model subfactor $R_u\subset R$ of index $4$. Then, $\{\frac{1}{4}\beta_m\,:\,m\in\bbn\}$ is a subset of $\Lambda(R_u\subset R)$, and all elements of this subset are arising from the finite-dimensional grid for the subfactor $R_u\subset R$.
\end{ppsn}
\begin{prf}
By \Cref{s}, we have $\beta_m=\Phi_4^{(2^m-1)}(1)=\frac{2^{m+1}-1}{2^m}\in\wedge_4\subseteq\Sigma_4$. By \Cref{known}, the matrix $\beta_m I_{4^r},\,r\in\bbn,$ can be written as sum of at most four projections in $M_{4^r}$ if and only if $\beta_m\in\Sigma_4$ and $4^r\beta_m\in\bbn$. Given $m\in\bbn$, choose $r_m\in\bbn$ such that $2r_m>m$. The condition $4^r\beta_m\in\bbn$ is then satisfied for $r=r_m$. Hence, $\beta_m I_{4^{r_m}}=p_1+p_2+p_3+p_4$ has a solution in the finite-dimensional Hilbert space $\bbc^{4^{r_m}}$. Therefore, $\frac{1}{4}\beta_m\in\Lambda(M_{4^{r_m}}\subset\Delta_4\otimes M_{4^{r_m}})\subset\Lambda(R_u\subset R)$, which concludes the proof.\qed
\end{prf}

Therefore, although for $R_u\subset R$ with $[R:R_u]=4$, the set $\Lambda(R_u\subset R)$ is completely known by Popa \cite{popa}, in \Cref{qp} we have identified a sequence of elements in $\Lambda(R_u\subset R)$ that arise from the finite-dimensional grid for $R_u\subset R$. With the help of this, we show the following result.

\begin{thm}\label{producing}
Consider a spin model subfactor $R_v\subset R$ such that $[R:R_v]=2n$, where $n\geq 3$, and let $\,t<\frac{1}{2}$ be given by the quadratic equation $\,t(1-t)=1/2n$. For $\alpha_m=\frac{1}{4}\Phi_4^{(2^m-1)}(1),\,m\in\bbn$, as in \Cref{qp}, we have $$\bigcup_{i=0}^{2n-4}\Big\{\frac{4\alpha_m+i}{2n}\,:\,m\in\bbn\Big\}\subseteq\Lambda(R_v\subset R)\cap(t,1-t)\,.$$ For each $0\leq i\leq 2n-4$, the sequence $\{(4\alpha_m+i)/2n\}_{m\geq 1}$ converges increasingly to $\frac{2+i}{2n}$.
\end{thm}
\begin{prf}
Consider any spin model subfactor $R_u\subset R$ such that $[R:R_u]=4$. Choose any $k\in\bbn$ and pick $\alpha\in\Lambda(M_{4^k}\subset \Delta_4\otimes M_{4^k})\subseteq\Lambda(R_u\subset R)$. By \Cref{qp}, we have $\alpha=\frac{1}{4}\Phi_4^{(2^m-1)}(1)$ for some $m\in\bbn$, and hence there are infinitely many choices for $\alpha$. Now, consider a spin model subfactor $R_v\subset R$ such that $[R:R_v]=2n$, where $n\geq 3$.
\medskip

\textbf{Step 1~:} $\frac{1}{2n}(4\alpha+i)$ lie in $\Lambda(R_v\subset R)$, where $0\leq i\leq 2n-4$ and $[R:R_v]=2n,\,n\geq 3$.
\smallskip

\noindent Since $\alpha\in\Lambda(M_{4^k}\subset \Delta_4\otimes M_{4^k})$, there are at most four projections $p_j,j=1,2,3,4,$ in $M_{4^k}$ such that $\sum_jp_j=4\alpha I_{4^k}$. By \Cref{known}, we have $4\alpha\in\Sigma_4$ and $4^{k+1}\alpha\in\bbn$. For $[R:R_v]=2n$, consider projections $q_j=I_{n^{2k}}\otimes p_j$ in $M_{(2n)^{2k}}$ for $j=1,2,3,4$. Then, $\sum_{j=1}^4q_j=4\alpha I_{(2n)^{2k}}$. Consider the following $2n-3$ projections in $\Delta_{2n}\otimes M_{(2n)^{2k}}\,:$
\[
\widetilde{q}_{0}=(q_1,q_2,q_3,q_4,0,\ldots,0)\,,
\]
and for $1\leq i\leq 2n-4$,
\[
\widetilde{q}_{i}=(q_1,q_2,q_3,q_4,I_{(2n)^{2k}},\ldots,I_{(2n)^{2k}},0,\ldots,0)\,.
\]
For the conditional expectation $E:\Delta_{2n}\otimes M_{(2n)^{2k}}\to M_{(2n)^{2k}}$, and for $0\leq i\leq 2n-4$, we get
\[
E(\widetilde{q}_{i})=\frac{1}{2n}\Big(\sum_{j=1}^4q_j+iI_{(2n)^{2k}}\Big)=\frac{1}{2n}(4\alpha+i)I_{(2n)^{2k}}\,.
\]
This says that $\frac{1}{2n}(4\alpha+i)$, for each $0\leq i\leq 2n-4$, are in the set $\Lambda(M_{(2n)^{2k}}\subset \Delta_{2n}\otimes M_{(2n)^{2k}})\subseteq\Lambda(R_v\subset R)$, where $[R:R_v]=2n$.
\medskip

\textbf{Step 2~:} $\frac{4\alpha}{6}>t$, where $t(1-t)=\frac{1}{6}$.
\smallskip

\noindent To show that $\frac{4\alpha}{6}>t$, we need to prove $\alpha>\frac{3}{4}\left(1-\frac{1}{\sqrt{3}}\right)$. Recall that $\alpha=\frac{1}{4}\Phi_4^{(2^m-1)}(1)$ for some $m\in\bbn$. By \Cref{s}, we have $\Phi_4^{(2^m-1)}(1)=\frac{2^{m+1}-1}{2^m}$. So, we need to show that $\frac{2^{m+1}-1}{2^m}>3-\sqrt{3}$ for any $m\in\bbn$. However, this follows easily since $2^m>\frac{1}{\sqrt{3}-1}$ for any $m\in\bbn$.
\medskip

\textbf{Step 3~:} $\frac{4\alpha+2}{6}<1-t$, where $t(1-t)=\frac{1}{6}$.
\smallskip

\noindent Observe that $t(1-t)=1/6$ and $t<\frac{1}{2}$ implies that $t=\frac{1}{2}\left(1-(1-\frac{4}{6})^{\frac{1}{2}}\right)$. Now, on the contrary, assume that $\frac{4\alpha+2}{6}\geq 1-t$. Then, we have the following\,:
\begin{IEEEeqnarray*}{lCl}
\frac{4\alpha+2}{6}\geq 1-t &\implies & \frac{4\alpha+2}{6}\geq\frac{1}{2}+\frac{1}{2\sqrt{3}}\\
&\implies & 4\alpha+2\geq 3+\sqrt{3}\\
&\implies & \alpha\geq\frac{1+\sqrt{3}}{4}\,.
\end{IEEEeqnarray*}
Since $\alpha\leq\frac{1}{2}$, we get that $\frac{1}{2}\geq\frac{1+\sqrt{3}}{4}$, that is, $2\geq 1+\sqrt{3}$, which is a contradiction. This proves that $\frac{4\alpha+2}{6}<1-t$.
\medskip

\textbf{Step 4~:} $\frac{4\alpha+i}{2n}$ lie in $(t,1-t)$, where $0\leq i\leq 2n-4$ and $t(1-t)=\frac{1}{2n},\,n\geq 3$.
\smallskip

\noindent It is enough to prove that $\frac{4\alpha}{2n}>t$ and $\frac{4\alpha+2n-4}{2n}<1-t$. We use induction on $n\geq 3$, where Step $3$ is the base step for the induction. Assume that the result is true for $2n-2$, and consider any spin model subfactor $R_v\subset R$ such that $[R:R_v]=2n$. First note that $t=\frac{1}{2}\left(1-(1-\frac{4}{2n})^{\frac{1}{2}}\right)$ as $t<1/2$ and $t(1-t)=1/2n$. Thus, we need to show that if
\begin{IEEEeqnarray}{lCl}\label{a}
\frac{4\alpha+2n-6}{2n-2}<\frac{1}{2}\left(1+\sqrt{1-\frac{4}{2n-2}}\right),
\end{IEEEeqnarray}
then
\begin{IEEEeqnarray}{lCl}\label{b}
\frac{4\alpha+2n-4}{2n}<\frac{1}{2}\left(1+\sqrt{1-\frac{4}{2n}}\right).
\end{IEEEeqnarray}
First we claim that the following inequality
\begin{IEEEeqnarray}{lCl}\label{d}
 \frac{n}{n-1}\left(1+\sqrt{1-\frac{2}{n}}\right)\geq 1+\sqrt{1-\frac{4}{2n-2}}+\frac{2}{n-1}=1+\sqrt{\frac{n-3}{n-1}}+\frac{2}{n-1}
\end{IEEEeqnarray}
holds. On the contrary, if we assume that
\[
\frac{n}{n-1}\left(1+\sqrt{1-\frac{2}{n}}\right)<1+\sqrt{\frac{n-3}{n-1}}+\frac{2}{n-1}\,,
\]
then we get the following inequality
\[
n+\sqrt{n(n-2)}<n+1+\sqrt{(n-1)(n-3)}\,.
\]
This further implies that
\[
n(n-2)+1-2\sqrt{n(n-2)}<(n-1)(n-3),
\]
which gives us $n-1<\sqrt{n(n-2)}$. However, this is clearly not possible for any $n\geq 3$. Therefore, \Cref{d} is validated. Now, to validate \Cref{b}, on the contrary, assume that $\frac{4\alpha+2n-4}{n}\geq 1+\sqrt{1-\frac{4}{2n}}$. Then, we have the following\,:
\begin{IEEEeqnarray*}{lCl}
& & \frac{4\alpha+2n-4}{n}\geq 1+\sqrt{1-\frac{4}{2n}}\\
&\implies& \frac{4\alpha+2n-6}{n-1}\geq \frac{n}{n-1}\left(1+\sqrt{1-\frac{2}{n}}\right)-\frac{2}{n-1}\\
&\implies& \frac{4\alpha+2n-6}{n-1}\geq 1+\sqrt{\frac{n-3}{n-1}}\qquad\mbox{by }\Cref{d}.
\end{IEEEeqnarray*}
In other words, we arrive at the inequality $\frac{4\alpha+2n-6}{n-1}\geq 1+\sqrt{1-\frac{4}{2n-2}}\,$, which contradicts to \Cref{a}. Therefore, \Cref{b} is validated. That is, $\frac{4\alpha+2n-4}{2n}<1-t$ for all $n\geq 3$, where $t(1-t)=1/2n$. Thus, the induction is complete, and we conclude that $\frac{4\alpha+i}{2n}<1-t$ for all $n\geq 3$ and $0\leq i\leq 2n-4$. The case of $\frac{4\alpha}{2n}>t$ is exactly similar, where Step $2$ is the base step for the induction, and we omit the details.
\smallskip

Finally, combination of Step $1$ and Step $4$ completes the first part of the proof. Now, for the second part, observe that the sequence $\{(4\alpha_m+i)/2n\}_{m\geq 1}$ is equal to $\{(2-\frac{1}{2^m}+i\big)/2n\}_{m\geq 1}$, which is clearly increasing. This follows directly from \Cref{s} by putting $k=2^m-1$. Indeed, $\Phi_4^{(2^m-1)}=\frac{2^{m+1}-1}{2^m}=2-\frac{1}{2^m}$ for $m\in\bbn$.\qed
\end{prf}
\smallskip

\noindent\textbf{Question G:} In the above theorem, we have produced a lot of rational numbers in $\Lambda(R_v\subset R)\cap (t,1-t)$ for the spin model subfactor $R_v\subset R$ with even index, but our method does not help to produce any irrational number. In fact, we do not know whether there is any irrational number in the set $\Lambda(R_v\subset R)\cap (t,1-t)$. We would like to propose this as a question, that is, whether the set $\Lambda(R_v\subset R)\cap (t,1-t)$ contain any irrational number for the cases of even index spin model subfactors.
\smallskip

Now, we turn our attention to the vertex model subfactors (recall from \Cref{preliminaries}). Recall that vertex model subfactors need not be irreducible, unlike the spin model subfactors. Below we produce infinitely many numbers in the set $\Lambda(R_w\subset R)\cap(t,1-t)$ for certain class of vertex model subfactors $R_w\subset R$. This is an important result of this paper because we move beyond `irreducibility'.

\begin{thm}\label{producing2}
Consider a vertex model subfactor $R_w\subset R$ such that $[R:R_w]=(2n)^2$, where $n\geq 3$, and let $\,t<\frac{1}{2}$ be given by the quadratic equation $t(1-t)=1/4n^2$. For $\alpha_m=\frac{1}{4}\Phi_4^{(2^m-1)}(1),\,m\in\bbn$, as in \Cref{qp}, we have
\[
\bigcup_{i=0}^{2n-4}\Big\{\frac{4\alpha_m+i}{2n}\,:\,m\in\bbn\Big\}\subseteq\Lambda(R_w\subset R)\cap(t,1-t)\,.
\]
For each $0\leq i\leq 2n-4$, the sequence $\{(4\alpha_m+i)/2n\}_{m\geq 1}$ converges increasingly to $\frac{2+i}{2n}$.
\end{thm}
\begin{prf}
Recall from \Cref{preliminaries}, the grid of finite-dimensional algebras for any vertex model subfactor $R_w\subset R$ of index $r^2$ is given by the following\,:
\[
\begin{matrix}
M_r\otimes\bbc &\subset & M_r\otimes M_r &\subset &\ldots &\subset & M_r\otimes^{\,k}M_r &\subset &\ldots &\subset  R\\
\cup & & \cup & & & &\cup & & &  \quad\cup\\
\bbc &\subset &\mbox{Ad}_w(\bbc\otimes M_r) &\subset &\ldots &\subset &\mbox{Ad}_{w^k}(\bbc\otimes^{\,k}M_r) &\subset &\ldots &\subset  R_w\\
\end{matrix}
\]
where $w$ is a bi-unitary in $M_{r^2}=M_r\otimes M_r$. The conditional expectation $E_k:M_r\otimes^{\,k}M_r\to\bbc\otimes^{\,k}M_r$ is given by $tr\otimes\mbox{id}$, where $tr$ is the unique normalized trace on $M_r$. Therefore, we have the following
\begin{IEEEeqnarray*}{lCl}
& & p_1+\ldots+p_r=r\alpha I_{r^k}\\
&\implies& E_k\big(E_{11}\otimes p_1+\ldots+E_{rr}\otimes p_r\big)=\alpha I_{r^k}
\end{IEEEeqnarray*}
for projections $p_j\in M_{r^k}$. This immediately says the following\,:
\begin{IEEEeqnarray}{lCl}\label{1}
\Lambda\big(\bbc\otimes^{\,k}M_r\subset \Delta_r\otimes^{\,k}M_r\big) &\subseteq& \Lambda\big(\bbc\otimes^{\,k}M_r\subset M_r\otimes^{\,k}M_r\big)
\end{IEEEeqnarray}
for any $k\in\bbn$ and fixed $r\in\bbn$. Now, let $t=\frac{1}{2}\left(1-(1-\frac{4}{r})^{\frac{1}{2}}\right)$ and $\widetilde{t}=\frac{1}{2}\left(1-(1-\frac{4}{r^2})^{\frac{1}{2}}\right)$. Then, clearly we have
\begin{IEEEeqnarray}{lCl}\label{2}
\widetilde{t}<t\qquad\mbox{and}\qquad 1-t<1-\widetilde{t}.
\end{IEEEeqnarray}
The result now follows immediately from \Cref{producing}, together with \Cref{1,2}, by putting $r=2n$.
\qed
\end{prf}

For any $m\in\bbn$ and each $0\leq i\leq 2n-4$, denote $\gamma_{m,i}:=\frac{4\alpha_m+i}{2n}$. For each $\gamma_{m,i}$, there is a sequence $\{\zeta_{m,i}^{(k)}\}_{k\geq 1}$ due to Proposition $5.5$ in \cite{popa}, defined recursively by $\zeta_{m,i}^{(0)}:=\gamma_{m,i}$ and $\zeta_{m,i}^{(k)}:=\Big(2n-\frac{1}{1-\zeta_{m,i}^{(k-2)}}\Big)^{-1}$. This is possible because if $N\subset M$ is either a spin model or a vertex model subfactor, then the inclusion $N\subset M$ splits $R$. Therefore, by \Cref{producing} and \ref{producing2}, for the cases of spin model subfactors of index $2n,n\geq 3,$ and vertex model subfactors of index $(2n)^2,n\geq 3$, we have the following infinite matrix where each entry is an element in the desired set $\Lambda(R_v\subset R)\cap (t,1-t)$,
\[
\begin{matrix}
\zeta_{1,i}^{(0)}=\gamma_{1,i} & & \zeta_{1,i}^{(1)} &\zeta_{1,i}^{(2)} &\cdots & \zeta_{1,i}^{(k)} &\cdots\\
 & & & & & \\
\zeta_{2,i}^{(0)}=\gamma_{2,i} & & \zeta_{2,i}^{(1)} &\zeta_{2,i}^{(2)} &\cdots & \zeta_{2,i}^{(k)} &\cdots\\
 & & & & & \\
\vdots &  &\vdots &\vdots & & \vdots\\
 & & & & & \\
\zeta_{m,i}^{(0)}=\gamma_{m,i} & & \zeta_{m,i}^{(1)} &\zeta_{m,i}^{(2)} &\cdots & \zeta_{m,i}^{(k)} &\cdots\\
 & & & & & \\
\vdots & &\vdots &\vdots & & \vdots\\
\end{matrix}
\]
for any $0\leq i\leq 2n-4$.
\smallskip

In the above infinite matrix, the elements $\gamma$ in the first column are constructed in \Cref{producing} and \ref{producing2}, and all the elements $\zeta$ in rest of the columns are produced from the first column recursively using Popa's construction (Proposition $5.5$, \cite{popa}) discussed in \Cref{Sec 3}. Moreover, $\lim_{m\to\infty}\gamma_{m,i}=\frac{2+i}{2n}$ for each $0\leq i\leq 2n-4$, and $\lim_{k\to\infty}\zeta_{m,i}^{(k)}=t$ for each $m\in\bbn$ and $0\leq i\leq 2n-4$. We conclude this section by showing that all the elements in the above infinite matrix are distinct.

\begin{ppsn}\label{popa sequence}
Fix any $0\leq i\leq 2n-4$, where $n\geq 3$. Then, we have $\zeta_{m_1,i}^{(k_1)}\neq\zeta_{m_2,i}^{(k_2)}$ for any $m_1\neq m_2\in\bbn$ or $k_1\neq k_2\in\bbn\cup\{0\}$.
\end{ppsn}
\begin{prf}
First observe that any two elements in a row are distinct due to the construction by Popa, that is, for each $m\in\bbn$, we have $\zeta_{m,i}^{(k_1)}\neq\zeta_{m,i}^{(k_2)}$ if $k_1\neq k_2$. Also, by our construction in \Cref{producing}, any two elements in the first column are distinct, that is, $\gamma_{m_1,i}\neq\gamma_{m_2,i}$ if $m_1\neq m_2$. Therefore, it is enough to show that $\zeta_{m_1,i}^{(k_1)}\neq\zeta_{m_2,i}^{(k_2)}$ for any $m_1\neq m_2$ and $k_1,k_2\in\bbn\cup\{0\}$.

First consider the case of $k_1=k_2=k$, that is, consider any two elements in the $k$-th column. Then, $\zeta_{m_1,i}^{(k)}=\zeta_{m_2,i}^{(k)}$ implies that $\gamma_{m_1,i}=\gamma_{m_2,i}$ due to their construction. Since this forces $m_1=m_2$, we conclude that any two elements in the $k$-th column, for $k\in\bbn$, are distinct. Now, assume that $\zeta_{m_1,i}^{(k_1)}=\zeta_{m_2,i}^{(k_2)}$ such that $m_1\neq m_2$ and $k_1\neq k_2$. Without loss of generality, we can take $k_1<k_2$. Again due to the construction, it readily follows that $\gamma_{m_1,i}=\zeta_{m_2,i}^{(k_2-k_1)}$. This says that it is enough to show that no element in the first column is equal to any element in any other column. Now taking $k_2-k_1=r$, we see that $\gamma_{m_1,i}=\zeta_{m_2,i}^{(r)}$ implies the following\,:
\[
2-\frac{1}{2^{m_1}}+i=\frac{2n}{2n-\frac{1}{1-\zeta_{m_2,i}^{r-2}}}\,.
\]
Therefore,
\[
\frac{2n}{2-\frac{1}{2^{m_1}}+i}=2n-\frac{1}{1-\zeta_{m_2,i}^{(r-2)}}\,,
\]
which further implies that
\[
\frac{1}{1-\zeta_{m_2,i}^{(r-2)}}=\frac{2n-\frac{2n}{2^{m_1}}+2ni}{2-\frac{1}{2^{m_1}}+i}\,.
\]
Therefore, we have the following equality
\[
1-\zeta_{m_2,i}^{(r-2)}=\frac{2^{m_1+1}-1+2^{m_1}i}{2n(2^{m_1}+1+2^{m_1}i)}\,.
\]
Since $\zeta_{m_2,i}^{(r-2)}<1/2$, we have the following inequality
\[
\frac{1}{2}<1-\zeta_{m_2,i}^{(r-2)}=\frac{2^{m_1+1}-1+2^{m_1}i}{2n(2^{m_1}+1+2^{m_1}i)},
\]
which further gives us the following\,:
\[
n<\frac{2^{m_1+1}-1+2^{m_1}i}{2^{m_1}+1+2^{m_1}i}\,.
\]
Therefore,
\[
2^{m_1}(n-2)+(n-1)2^{m_1}i)<-1-n
\]
which is impossible since $n\geq 3$. Hence, $\zeta_{m_1,i}^{(k_1)}\neq\zeta_{m_2,i}^{(k_2)}$ for $m_1\neq m_2$ and $k_1\neq k_2$, which concludes the proof.\qed
\end{prf}

Combining Theorems \ref{aa}, \ref{bb}, Theorems \ref{producing}, \ref{producing2}, and \Cref{popa sequence}, we conclude the following final result.
\begin{thm}
\begin{enumerate}[$(i)$]
\item Let $R_v\subset R$ be a spin model subfactor of index $2n,\,n\geq 3$. The following is a subset of the set $\Lambda(R_v\subset R)\cap(t,1-t)$, where $\,t(1-t)=1/2n$,
\begin{IEEEeqnarray*}{lCl}
& & \{m/2n\,:\,0\leq m\leq 2n\}\cup\{1/k\,:\,2\leq k\leq 2n-2\}\cup\{(k-1)/k\,:\,2\leq k\leq 2n-2\}\\
& & \cup\,\{\gamma_{m,i}\,:m\in\bbn,\,0\leq i\leq 2n-4\}\cup\big\{\zeta^{(k)}_{m,i}\,:m,k\in\bbn,\,0\leq i\leq 2n-4\big\}\\
& & \cup\,\{1-\gamma_{m,i}\,:m\in\bbn,\,0\leq i\leq 2n-4\}\cup\big\{1-\zeta^{(k)}_{m,i}\,:m,k\in\bbn,\,0\leq i\leq 2n-4\big\}.
\end{IEEEeqnarray*}
\item Let $R_w\subset R$ be a vertex model subfactor of index $(2n)^2,\,n\geq 3$. The following is a subset of the set $\Lambda(R_w\subset R)\cap(t,1-t)$, where $t(1-t)=1/(2n)^2$,
\begin{IEEEeqnarray*}{lCl}
& & \{m/4n^2\,:\,0\leq m\leq 4n^2\}\cup\{1/4n^2k\,:\,1\leq k\leq 2n-1\}\cup\{(4n^2k-1)/k\,:\,1\leq k\leq 2n-1\}\\
& & \cup\,\{\gamma_{m,i}\,:m\in\bbn,\,0\leq i\leq 2n-4\}\cup\big\{\zeta^{(k)}_{m,i}\,:m,k\in\bbn,\,0\leq i\leq 2n-4\big\}\\
& & \cup\,\{1-\gamma_{m,i}\,:m\in\bbn,\,0\leq i\leq 2n-4\}\cup\big\{1-\zeta^{(k)}_{m,i}\,:m,k\in\bbn,\,0\leq i\leq 2n-4\big\}\,.
\end{IEEEeqnarray*}
\end{enumerate}
\end{thm}

%%%%%%%%%%%%%%%%%%%%%%%%%%%%%%%%%%%%%%%%%%
%%%%%%%%%%%%%%%%%%%%%%%%%%%%%%%%%%%%%%%%%%
%%%%%%%%%%%%%%%%%%%%%%%%%%%%%%%%%%%%%%%%%%\

\section*{Acknowledgements}
Keshab Chandra Bakshi acknowledges the support of DST INSPIRE Faculty fellowship\\ DST/INSPIRE/04/2019/002754, and Satyajit Guin acknowledges the support of SERB grant
MTR/2021/000818.

%%%%%%%%%%%%%%%%%%%%%%%%%%%%%%%%%%%%%%%%
%%%%%%%%%%%%%%%%%%%%%%%%%%%%%%%%%%%%%%%%
%%%%%%%%%%%%%%%%%%%%%%%%%%%%%%%%%%%%%%%%

\bigskip

\bigskip

\noindent{\sc Keshab Chandra Bakshi} (\texttt{keshab@iitk.ac.in, bakshi209@gmail.com})\\
         {\footnotesize Department of Mathematics and Statistics,\\
         Indian Institute of Technology, Kanpur,\\
         Uttar Pradesh 208016, India}
\bigskip

\noindent{\sc Satyajit Guin} (\texttt{sguin@iitk.ac.in})\\
         {\footnotesize Department of Mathematics and Statistics,\\
         Indian Institute of Technology, Kanpur,\\
         Uttar Pradesh 208016, India}

\end{document}